\documentclass[12pt]{article}
\usepackage[cp866]{inputenc}
\usepackage{citehack}
\usepackage{amsmath}
\usepackage{amsfonts}
\usepackage{amssymb}
\newtheorem{theorem}{Theorem}
\newtheorem{prop}{Proposition}

\def\spa{\mathop\text{{\rm span}}\nolimits}

\def\Hom{\mathop\text{\rm Hom}\nolimits}

\def\tr{\mathop\text{\rm tr}\nolimits}

\def\Real{\mathbb{R}}

\def\g{\mathfrak{g}}
\def\h{\mathfrak{h}}
\def\hol{\mathfrak{hol}}
\def\so{\mathfrak{so}}

\def\gl{\mathfrak{gl}}

\def\f{\mathfrak{f}}

\def\z{\mathfrak{z}}
\def\R{\mathcal{R}}

\def\G{\Gamma}
\def\F{\mathcal F}
\def\T{T^\bot}
\def\n{\nabla}
\def\nb{\bar\nabla}
\def\nt{\nabla^t}
\def\ns{\overset{*}{\nabla}}
\def\hs{\overset{*}{h}}
\def\As{\overset{*}{A}}
\def\Rs{\overset{*}{R}}
\def\nst{\overset{*}{\nabla}^t}

\def\M{\bar M}
\def\id{\mathop\text{\rm id}\nolimits}
\def\pr{\mathop\text{\rm pr}\nolimits}

\title{Lightlike  foliations on Lorentzian manifolds with weakly irreducible holonomy algebra}

\author{Natalia Bezvitnaya}

\begin{document}

\maketitle

\maketitle\vskip-50ex
 {\renewcommand{\abstractname}{Abstract}\begin{abstract}
We study the lightlike foliations that appear on Lorentzian manifolds with weakly irreducible not
irreducible holonomy algebra. We give global structure equations for the foliation that generalize
the Gauss and Weingarten equations for one lightlike hypersurface. This gives us some global
operators on the manifold. Using these operators, we decompose the curvature tensor of the manifold into several components.
We give a criteria how to find the type of the holonomy algebras (there are 4 possible types) in terms of the global
operators.
\end{abstract}



\section*{Introduction}

Weakly irreducible not irreducible holonomy algebras of Lorentzian manifolds were classified recently. In \cite{B-I} L. Berard Bergery
and A. Ikemakhen divided these holonomy algebras into 4 types and associated to each such algebra a subalgebra of $ \so(n) $ that is 
called the orthogonal part (the dimension of the manifold  is $ n+2 $).  In \cite{LeD,Le3} T. Leistner proved that the orthogonal
part of the holonomy algebra of a Lorentzian manifold
is the holonomy algebra  of a Riemannian manifold. We recall this classification in section \ref{hol}.

If the holonomy algebra of a connected Lorentzian manifold $\bar{M}$ is weakly irreducible and not irreducible (i.e. it preserves an 
isotropic line  and does not preserve any nondegenerate vector subspace of the tangent space), then we obtain  on $\bar{M}$ a parallel 
distribution $D$ of isotropic lines and the distribution $D^\bot$
 of degenerate hypersubspaces of the tangent spaces. The distribution 
$D^\bot$ is parallel and,  in particular,  it is involutive. This gives us a foliation in lightlike hypersurfaces.

Our purpose is to describe the  geometry of Lorentzian manifolds with weakly irreducible not irreducible holonomy algebras of each
type  and with each orthogonal part in terms of the foliation in lightlike hypersurfaces.
There were some approaches to study geometry of these manifolds in local coordinates, see \cite{Wal,Sch,Ik,LeD,Bou}. 
We hope that our approach can give some global description.

The case of one lightlike hypersurface of a Lorentzian manifold was studied, for example, in \cite{DB}. In section \ref{S1} we recall
some results from this book. Note that in order to obtain a connection on a lightlike hypersurface $M$ of a Lorentzian manifold
$\bar M$ we must choose a distribution (a screen distribution) $S(TM) \subset TM$  
that is any complement distribution to $TM^\bot\subset TM$ or, 
equivalently, choose a complementary isotropic vector bundle (a transversal vector bundle) $tr(TM)$ to $TM\subset T\bar{M}|_M$.
Such choice is not always canonical. Some approaches to a canonical choice of $S(TM)$ (equivalently of $tr(TM)$) can be found 
in \cite{AG}.

In the beginning of section \ref{S3} we rewrite the structure equations from \cite{DB} for one hypersurface of the foliation.
Since we have not one lightlike hypersurface but a foliation in lightlike hypersurfaces, it is natural  to have global structure
equations. Such equations can be obtained by choosing  a global  screen distribution (or transversal bundle), the existence
of the last is guaranteed  by theorem \ref{Tr}.
From the global equations we obtain some operators on $\bar{M}$ (which depend on the choice of the screen distribution). These 
operators generalize the operators for one hypersurface from \cite{DB}.

A global  screen distribution  gives us a decomposition of the tangent space of $\bar{M}$ at each point. The same decomposition 
was used in \cite{G} in order to decompose the curvature tensor  at a point into several  components, see section \ref{G}.
Knowing these components at each point of the manifold, we know the type of the holonomy algebra.
In section \ref{R} we consider these components for the curvature tensor of the manifold $\M$
 and express them in terms of the global operators that we have defined. 

In section \ref{type}   we give a criterion  how to find the type of the holonomy algebra in terms of our global operators. 
First we  distinguish  holonomy algebras of type 2 and 4 from the holonomy algebras of type 1 and 3.
Then we give criterions for holonomy algebras of type 3 and 4.

Finally we consider a local example, where the holonomy algebra is Abelian.

Another open problem is to find a canonical way of choosing the screen distribution.        

{\it Acknowledgments} I would like to thank Helga Baum for the theme proposed to me and for
the useful discussions during October-December 2004.

\section{Case of one lightlike hypersurface of\\ a Lorentzian manifold}\label{S1}

In this section we recall some results from \cite{DB}.

Let $(\bar{M},\bar{g})$ be a connected  Lorentzian manifold of dimension $n+2$ and
let $M\subset \bar M$ be a {\it lightlike} hypersurface of $\bar M$, i.e.
the restriction of $\bar g_x$ to $T_xM$ is degenerate for all $x\in M$.
Since $T_xM$ is degenerate, the perpendicular $T_xM^\bot$ to $T_xM$ in $T_x\bar M$
is an isotropic line, which  is contained in $T_xM$. We get on $M$ the vector subbundle $TM^\bot=\cup_{x\in M}T_xM^\bot$ of $TM$.

In order to obtain a connection on $M$ we fix a complementary
vector bundle $S(TM)$ of $TM^\bot$ in $TM$. We have \begin{equation}\label{l1}  TM=S(TM)\overset{\bot}{\oplus}TM^\bot. \end{equation}
If $M$ is paracompact, $S(TM)$ always exists. The distribution  $S(TM)\subset TM$ is called a {\it screen distribution}.  

Since any maximal isotropic subspace of $T_x\bar M$ is 1-dimensional, we see that the distribution $S(TM)$ is not degenerate.
Hence we obtain the decomposition \begin{equation}\label{l2}  T\bar M|_M=S(TM)\overset{\bot}{\oplus}S(TM)^\bot, \end{equation}
where $S(TM)^\bot$ is the orthogonal complementary vector  bundle to $S(TM)$ in $T\bar M_{|M}$.

In \cite{DB} it was proved that for a given screen distribution $S(TM)$  there exists a unique vector bundle $\tr(TM)$ of rank $1$ 
over $M$, such that for any non-zero section
$\xi$ of $TM^\bot$ on a coordinate neighborhood $U\subset M$, there exists a unique section $N$ of $\tr(TM)$ on $U$ such that
$$\bar g(\xi,N)=1 \text{ \rm  and  } \bar g(N,N)=\bar g(N,W)=0 \text{ \rm   for all } W\in \Gamma(S(TM)_{|U}).$$
The vector bundle $\tr(TM)$ is called the {\it lightlike transversal vector bundle} of $M$ with respect to $S(TM)$.  
We have the following decompositions of $T\bar M_{|M}$:
\begin{equation}\label{l3}  T\bar M_{|M}=S(TM)\overset{\bot}{\oplus}(TM^\bot\oplus \tr(TM))=TM\oplus\tr(TM). \end{equation}

Suppose that we have a screen distribution on a lightlike hypersurface $M$ of a Lorentzian manifold $(\M,\bar g)$.
Using the second form of the decomposition \eqref{l3}, we obtain
\begin{equation}\label{l4} \nb_XY=\underbrace{\n_XY}_{\in \G(TM)}+\underbrace{h(X,Y)}_{\in \G(\tr(TM))} \end{equation}  
and
\begin{equation}\label{l5} \nb_XV=\underbrace{-A_VX}_{\in \G(TM)}+\underbrace{\nt_XV}_{\in \G(\tr(TM))} \end{equation}
for any $X,Y\in\G(TM)$ and $V\in\G(\tr(TM)).$ It is known that $\n$ is a torsion-free linear connection on $M$,
$h$ is a $\G(\tr(TM))$-valued symmetric $\F(M)$-bilinear form on $\G(TM)$, $A_V:\G(TM)\to \G(TM)$ is a $\F(M)$-linear operator and
$\nt$ is a linear connection on the vector bundle $\tr(TM)$.    
The connections $\n$ and $\nt$ are called {\it the induced connections on} $M$ and $\tr(TM)$ respectively. The bilinear form $h$ and the operator
$A_V$ are called {\it the second fundamental form} and {\it the shape operator} respectively. The formulas \eqref{l4} and \eqref{l5}
are called {\it the Gauss} and {\it Weingarten formulas} respectively.

Using decomposition \eqref{l1}, we obtain
\begin{equation}\label{l6} \n_XY=\underbrace{\ns_XY}_{\in \G(S(TM))}+\underbrace{\hs(X,Y)}_{\in \G(TM^\bot)} \end{equation}  
and
\begin{equation}\label{l7} \n_XU=\underbrace{-\As_UX}_{\in \G(S(TM))}+\underbrace{\nst_XU}_{\in \G(TM^\bot)} \end{equation}
for any $X\in\G(TM)$, $Y\in \G(S(TM))$ and $U\in\G(TM^\bot).$
We have  that $\ns$ and $\nst$  are linear connections on vector bundles $S(TM)$ and $TM^\bot$ respectively, 
$\hs:\G(TM)\times \G(S(TM))\to \G(TM^\bot)$ is a  $\F(M)$-bilinear form and  $\As_U:\G(TM)\to \G(S(TM))$ is a $\F(M)$-linear operator.
The bilinear form $\hs$ and the operator
$\As_U$ are called {\it the second fundamental form} and {\it the shape operator of the screen} distribution $S(TM)$ respectively.
 The equations \eqref{l6} and \eqref{l7}
are called {\it the Gauss} and {\it Weingarten equations for the screen distribution} respectively.

Let $\bar R$ and $R$ be the curvature tensors for the connections $\nb$ and $\n$ respectively.
Then for $X,Y,Z\in \G(TM)$ we have
\begin{equation} \label{l8} \bar R(X,Y)Z=R(X,Y)Z+A_{h(X,Z)}Y-A_{h(Y,Z)}X+(\n_Xh)(Y,Z)-(\n_Yh)(X,Z), \end{equation}
where $$(\n_Xh)(Y,Z)=\nt_X(h(Y,Z))-h(\n_XY,Z)-h(Y,\n_XZ).$$
The equation \eqref{l8} is called {\it the Gauss-Codazzi equation.}

\section{Global structure equations}\label{S3}

Now consider a case when we obtain a foliation of lightlike hypersurfaces on a Lorentzian manifold. 

Suppose that the holonomy group of a connected $n+2$-dimensional  Lorentzian manifold $(\M,\bar g)$ is weakly irreducible 
and not irreducible, i.e. for any $x\in\M$ the holonomy group $Hol_x\subset O(T_x\M,\bar g_x)$ preserves an
isotropic line of $T_x\M$ and does not preserve any nondegenerate vector subspace of $T_x\M$. 
Then we obtain on $\M$ a parallel distribution $D$ of isotropic lines. Since $Hol_x\subset O(T_x\M,\bar g_x)$,
we see that the perpendicular $D^\bot_x\subset T_x\M$ is also preserved. Obviously, $D_x\subset D^\bot_x$,  $\dim D^\bot=n+1$
and the subspace $D^\bot$ is degenerate. We obtain a parallel distribution $D^\bot$ on $\M$. It is known that the distributions
$D$ and $D^\bot$ are smooth. Denote $D^\bot$ by $T$,  then $D=T^\bot$. 
Since the torsion of the Levi-Civita connection is zero, for any $X,Y\in\G(T)$ we have $[X,Y]=\n_XY-\n_YX\in\G(T)$,
hence the distribution $T$ is involutive.    
By the Frobenius theorem, for any point $x\in\M$ we have a maximal integral manifold $M_x$ of the distribution $T$.
Obviously, $M_x\subset\M$ is a lightlike hypersurface. Thus we obtain a foliations of lightlike hypersurfaces on $\M$.

Let $x\in \M$. We will rewrite the above formulas for the lightlike hypersurface $M_x\subset \M$. Suppose that we have some screen
distribution $S(TM_x)$. Since the distribution $T$ is parallel, from \eqref{l4} we obtain
\begin{equation}\label{l9} \nb_XY=\n_XY\in \G(TM_x) \end{equation}
for all  $X,Y\in \G(TM_x)$, i.e. $h=0$.

The equations \eqref{l5} and \eqref{l6} remain without change.   
From \eqref{l7} we obtain
\begin{equation}\label{l10} \n_XU=\nst_XU\in \G(TM^\bot_x) \end{equation}
for all  $X\in\G(TM_x)$, $U\in\G(TM^\bot_x)$, i.e. $\As=0$.

Since $h=0$, from \eqref{l8} it follows that
for $X,Y,Z\in \G(TM_x)$ we have
\begin{equation} \label{l11} \bar R(X,Y)Z=R(X,Y)Z. \end{equation}
 
Since we have a foliation of lightlike hypersurfaces, it is natural to construct a global screen distribution
and the corresponding transversal bundle. 


\begin{theorem}\label{Tr} For the distribution $\T$ there  exists   on $\M$ a (not unique) distribution $Tr$ of rank 1
that satisfies the condition: for any  $\xi\in\G(\T)$ nonzero at all points of some open set $U\subset\M $, there exists a unique section $N\in\G(Tr_{|U})$
satisfying $\bar g(N,N)=0$ and  $\bar g(\xi,N)=1$ on $U$.
\end{theorem}  
{\it Proof.} 
In \cite{Sch} it was proved that locally there exist coordinates $x^0,x^1,...x^n,x^{n+1}$  
such that the metric $\bar g$ has the form $$\bar g=2dx^0dx^{n+1}+\sum_{i,j=1}^{n}g_{ij}dx^idx^j+f\cdot(dx^{n+1})^2,$$
where $g_{ij}$ are functions of $x^1,...,x^n,x^{n+1}$ and $f$ is a function of $x^0,...,x^n,x^{n+1}$ .
For all $x$ from this coordinate neighborhood we have $\T_x=\Real(\frac{\partial}{\partial x^{0}})_x$ and 
we choose $Tr_x=\Real(-\frac{1}{2}f(\frac{\partial}{\partial x^0})_x+(\frac{\partial}{\partial x^{n+1}})_x)$. Thus  on some neighborhood
of each point $x\in\M$  we obtain a distribution that satisfy the theorem. We assume that $\M$ 
is paracompact, then we have a locally finite open covering $(U_i)_{i\in I}$ of $\M$, where
$U_i$ are coordinate neighborhoods as above. On each $U_i$ we have a distribution $Tr_i$
that satisfy the theorem. Let $(f_i)_{i\in I}$ be a partition of unity for the covering $(U_i)_{i\in I}$.
Let $x\in\M$. Since $(U_i)_{i\in I}$ is locally finite, the point $x$ is contained in a finite number of the open sets
$U_{i_1},...,U_{i_m}$, where $i_1,...,i_m\in I$. Let $V=U_{i_1}\cap...\cap U_{i_m}$.
Let $\xi\in\G(\T)$ be any  section nonzero at all points of $V$. Then for each $1\leq k\leq m$ there exists
a unique $N_{i_k}\in\G({Tr_{i_k}}_{|U_{i_k}})$ with $\bar g(\xi,N_{i_k})=1$ on $V$.
For $N=f_{i_1}N_{i_1}+\cdots+f_{i_m}N_{i_m}-h\xi$, where $2h=\bar g(f_{i_1}N_{i_1}+\cdots+f_{i_m}N_{i_m},f_{i_1}N_{i_1}+\cdots+f_{i_m}N_{i_m})$,
 we have $\bar g(N,N)=0$ and
 $$\bar g(\xi,N)=\bar g(\xi,f_{i_1}N_{i_1}+\cdots+f_{i_m}N_{i_m})=f_{i_1}\bar g(\xi,N_{i_1})+\cdots+f_{i_m}\bar g(\xi,N_{i_m})=
f_{i_1}+\cdots+f_{i_m}=1.$$
For any $y\in V$ we define $Tr_y$ as $Tr_y=\Real N_y$. We do this for all $x\in\M$ and obtain a distribution $Tr$ on $\M$. 
Obviously, $Tr$ is well defined  on $\M$ and satisfies the theorem.
$\Box$

We see that the vector subspace $\T_x\oplus Tr_x\subset T_x\M$ is nondegenerate, hence the orthogonal complement
$S_x=(\T_x\oplus Tr_x)^\bot\subset T_x\M$ is an Euclidean subspace. Hence the restriction of $\bar g$ on
the distribution $S=(\T\oplus Tr)^\bot\subset T\M$ is positively definite. We have
\begin{equation}\label{l12}  T=S\overset{\bot}{\oplus}T^\bot \end{equation}
and \begin{equation}\label{l13}  T\M=S\overset{\bot}{\oplus}(T^\bot\oplus Tr)=T\oplus Tr. \end{equation}

Let $M_x\subset \M$ be an integral manifold of $T$ through $x\in\M$, then $S(TM_x)=S_{|M_x}$ is a
screen distribution on $M_x$ and $\tr(TM_x)=Tr_{|M_x}$ is the corresponding transversal bundle on $M_x$.

Now we can generalize the operators $\n$, $A_V$, $\nt$, $\hs$, $\ns$ and $\nst$.

Since the distribution $\T$ is parallel, we have 
\begin{equation}\label{l14} \nb_WY=\n_WY\in \G(T) \end{equation}  
Using the second form of the decomposition \eqref{l13}, we obtain 
\begin{equation}\label{l15} \nb_WV=\underbrace{-A_VW}_{\in \G(T)}+\underbrace{\nt_WV}_{\in \G(Tr)} \end{equation}
for any $W\in\G(T\M)$, $Y\in\G(T)$ and $V\in\G(Tr).$ It can be proved that $\n$ is a linear connection on the bundle $T$,
$A_V:\G(T\M)\to \G(T)$ is a $\F(\M)$-linear operator, and
$\nt$ is a linear connection on the vector bundle $Tr$.    
We call $\n$ and $\nt$ {\it the induced connections on} $T$ and $Tr$ respectively.
We call $A_V$ {\it the shape operator} of $\M$ with respect to the screen distribution $S$.

Using decomposition \eqref{l12}, we obtain
\begin{equation}\label{l16} \n_WY=\underbrace{\ns_WY}_{\in \G(S)}+\underbrace{\hs(W,Y)}_{\in \G(T^\bot)} \end{equation}  
and
\begin{equation}\label{l17} \n_WU=\nst_WU\in \G(T^\bot) \end{equation}
for any $W\in\G(T\M)$, $Y\in \G(S)$ and $U\in\G(\T).$
It can be proved that   $\ns$ and $\nst$  are linear connections on vector bundles $S$ and $T^\bot$ respectively, 
$\hs:\G(T\M)\times \G(S)\to \G(T^\bot)$ is a  $\F(M)$-bilinear form. We call $\hs$ {\it the second fundamental form} of
the screen distribution $S$.

For a hypersurface $M_x$ the operators $\n$, $A_V$, $\nt$, $\hs$, $\ns$ and $\nst$ can be obtained from the 
defined above  by taking the obvious restrictions. Restricting to $M_x$ \eqref{l14}, \eqref{l15}, \eqref{l16} and \eqref{l17}
we obtain \eqref{l9}, \eqref{l5}, \eqref{l6} and \eqref{l10}  respectively. 

Let $V\in \G(Tr)$. We have $\bar g(V,V)=0$, hence for any $W\in \G(T\M)$ holds $\bar g(\nb_WV,V)+\bar g(V,\nb_WV)=0$ and
\begin{equation} \bar g(\nb_WV,V)=0. \end{equation}
From \eqref{l15} we obtain $\bar g(-A_VW+\nt_WV,V)=0$. Thus, \begin{equation}\bar g(A_VW,V)=0.\end{equation}
This means, that $A_V$ takes values in $\G(S)$. Thus we can consider $A_V$ as
\begin{equation} A_V:\G(T\M)\to\G(S).\end{equation}

By analogy, for any $Y\in\G(S)$ and $V\in\G(Tr)$ we have $\bar g(Y,V)=0$, hence for any $W\in \G(T\M)$ holds $\bar g(\nb_WY,V)+\bar g(Y,\nb_WV)=0$.
Using \eqref{l15} and \eqref{l16}, we obtain
\begin{equation}\label{l18} \bar g(\hs(W,Y),V)=\bar g(A_VW,Y) \end{equation} 
for all $W\in\G(T\M)$,  $Y\in\G(S)$ and $V\in\G(Tr)$.
Hence $A_V$ and $\hs$ can be found from each other.

\section{Weakly irreducible not irreducible holonomy algebras of Lorentzian manifolds}\label{hol}

Let $(\Real^{1,n+1},\eta)$ be a Minkowski space of dimension $n+2$,
where  $\eta$ is a metric on $\Real^{n+2}$ of
signature $(1,n+1)$. We fix a basis
$U,X_1,...,X_n,V$ of $\Real^{1,n+1}$ with respect to which the Gram matrix
of $\eta$ has the form
$\left(\begin{array}{ccc}
0 & 0 & 1\\ 0 & E_n & 0 \\ 1 & 0 & 0 \\
\end{array}\right)$,
 where $E_n$ is the $n$-dimensional identity matrix.

A subalgebra $\g\subset \so(1,n+1)$ is
called {\it irreducible} if it does not preserve any proper
subspace of $\Real^{1,n+1}$; $\g$ is called {\it weakly irreducible} if
it does not  preserve any nondegenerate proper subspace of $\Real^{1,n+1}$.
Obviously, if $\g\subset \so(1,n+1)$ is
irreducible, then it is weakly irreducible.
If $\g\subset \so(1,n+1)$ preserves
a degenerate proper subspace $W\subset \Real^{1,n+1}$, then it preserves the
isotropic line $W\cap W^\bot$. 

From the classification of M. Berger (see \cite{Ber}) it follows that the 
only irreducible holonomy algebra of Lorentzian manifolds is isomorphic to $\so(1,n+1)$,
so we consider only the case of weakly irreducible not irreducible holonomy algebra.

Let $(\M,\bar g)$ be an $n+2$-dimensional connected Lorentzian manifold with weakly irreducible not
irreducible holonomy algebra. Let $x\in\M$. We identify $(T_x\M,\bar g_x)$ with  $(\Real^{1,n+1},\eta)$.   
We assume that the subalgebra $\hol\subset\so(1,n+1)$ corresponding to the holonomy algebra
$\hol_x\subset\so(T_x\M,\bar g_x)$ preserves the isotropic line $\Real U$, i.e. $\hol$ is
contained in  the subalgebra $\so(1,n+1)_{\Real U}$  of $\so(1,n+1)$ that preserves the line $\Real U$.
Above we had a decomposition $T_x\M=\T_x\oplus S_x\oplus Tr_x$. We assume that $X_1,...,X_n$ correspond
to a basis of $S_x$ and $V$ corresponds to a vector of $Tr_x$.

The Lie algebra $\so(1,n+1)_{\Real U}$ can be identified with
the following algebra of matrices:
  $$\so(1,n+1)_{\Real U}=\left\{ \left (\begin{array}{ccc}
a &X & 0\\ 0 & A &-X^t \\ 0 & 0 & -a \\
\end{array}\right):\, a\in \Real,\, X\in \Real^n,\,A \in \so(n)
 \right\}.$$

Let $\h\subset\so(n)$ be a subalgebra. Recall that $\h$ is a compact Lie
algebra and we have the decomposition $\h=\h'\oplus\z(\h)$, where
$\h'$ is the commutant of $\h$ and $\z(\h)$ is the center of $\h$.

The following result is due to L. Berard Bergery and A. Ikemakhen.

{\bf Theorem}
{\it Suppose $\hol\subset \so(1,n+1)_{\Real U}$ is a weakly irreducible holonomy algebra.
Then $\hol$ belongs to one of the following types
\begin{description}
\item[type 1.] $\hol^{1,\h}=\left\{ \left (\begin{array}{ccc}
a &X & 0\\ 0 & A &-X^t \\ 0 & 0 & -a \\
\end{array}\right):\, a\in \Real,\, X\in \Real^n,\,A \in \h
 \right\}$, where
$\h\subset\so(n)$ is a subalgebra;

\item[type 2.] $\hol^{2,\h}=\left\{ \left (\begin{array}{ccc}
0 &X & 0\\ 0 & A &-X^t \\ 0 & 0 & 0 \\
\end{array}\right):\,  X\in \Real^n,\,A \in \h
 \right\}$;
\item[type 3.] $\hol^{3,\h,\varphi}=\left\{ \left (\begin{array}{ccc}
\varphi(B) &X & 0\\ 0 & A+B &-X^t \\ 0 & 0 & -\varphi(B) \\
\end{array}\right):\,  X\in \Real^n,\,A \in \h',\,B\in\z(\h) \right\}$,
where $\varphi :\h\to\Real$ is a non-zero linear map with $\varphi|_{\h'}=0$;

\item[type 4.] $\hol^{4,\h,\psi}=\left\{ \left (\begin{array}{cccc}
0 &X&\psi(B) & 0\\ 0 & A+B&0 &-X^t \\ 0 & 0 & 0 &-\psi(B)^t \\
0&0&0&0\\
\end{array}\right):\,  X\in \Real^{n_1},\,A \in \h',B\in\z(\h) \right\}$,
where we have a non-trivial   
decomposition $n=n_1+n_2$   such that $\h\subset\so(n_1)$;
 and $\psi:\h\to \Real^{n_2}$ is a surjective linear map with $\psi|_{\h'}=0$.
\end{description}}

The subalgebra $\h\subset\so(n)$ associated to a holonomy algebra in the above theorem is called {\it the
orthogonal part} of the holonomy algebra. 
In \cite{Le3,LeD} T. Leistner proved the following theorem 

{\bf Theorem} {\it The orthogonal part of the weakly irreducible not irreducible holonomy algebra of a Lorentzian
manifold is the holonomy algebra of a Riemannian manifold.}

\section{Decomposition of an abstract curvature tensor} \label{G}

In this section we use the notations of section \ref{hol}. We recall the decomposition of a curvature tensor
for the holonomy algebras given in \cite{G}. 

Let $W$ be a vector space and $\f\subset \gl(W)$ a subalgebra. The space of curvature tensors for the Lie  algebra $\f$
is defined as follows
$$\R(\f)=\{R\in\text{Hom}(W\otimes W,\f):R(u,v)=-R(v,u),
R(u,v)w+R(v,w)u+R(w,u)v=0 \text{ }\text{\it for all } u,v,w\in W\}.$$

Suppose that the holonomy algebra of $\M$ is of type 1, i.e. $\hol=\hol^{1,\h}$ for some $\h\subset\so(n)$.
Let $\bar R_x$ be the curvature tensor of the manifold $\M$ at a point $x\in\M$, then  $\bar R_x\in\R(\hol^{1,\h})$, where $x\in \M$
(we use the identifications as in section \ref{hol}).
Having the decomposition $\Real^{1,n+1}=\Real U\oplus\Real^n\oplus\Real V$, $\bar R_x$ can be decomposed
into 5 components,  $\bar R_{x}=\bar R_{x1}+\bar R_{x2}+\bar R_{x3}+\bar R_{x4}+\bar R_{x5}$.
These components are given by  elements
$$R_\h\in\R(\h), \text{  } P\in\{P\in \Hom (\Real^n,\h):
\,\eta(P(u)v,w)+\eta(P(v)w,u)+\eta(P(w)u,v)=0\text{ }\text{\it for all }u,v,w\in \Real^n\},$$
$$T\in \Hom(\Real^n,\Real^n) \text{ with } T^*=T, \text{  } L\in\Hom(\Real^n,\Real) \text{  and } \lambda\in\Real,$$
and can be found from the following  conditions
$$
\begin{array}{llll}
\bar R_{x1}(X,Y)&=  \left (\begin{array}{ccc}
0 &0 & 0\\ 0 & R_h(X,Y) &0 \\ 0 & 0 & 0 \\
\end{array}\right),& & \\
\bar  R_{x2}(V,X)&=  \left (\begin{array}{ccc}
0 &0 & 0\\ 0 & P(X) &0 \\ 0 & 0 & 0 \\
\end{array}\right), &
\bar R_{x2}(X,Y)&=  \left (\begin{array}{ccc}
0 &P^*(X,Y) & 0\\ 0 & 0 &-P^*(X,Y)^t \\ 0 & 0 & 0 \\
\end{array}\right),\\ 
\bar R_{x3}(V,X)&=  \left (\begin{array}{ccc}
0 &T(X) & 0\\ 0 & 0 &-T(X)^t \\ 0 & 0 & 0 \\
\end{array}\right),\\ 
\bar R_{x4}(U,V)&=  \left (\begin{array}{ccc}
\lambda &0 & 0\\ 0 & 0 &0 \\ 0 & 0 & -\lambda \\
\end{array}\right),& & \\
\bar R_{x5}(V,X)&=  \left (\begin{array}{ccc}
L(X) &0 & 0\\ 0 & 0 &0 \\ 0 & 0 & -L(X) \\
\end{array}\right),& 
\bar R_{x5}(U,V)&=  \left (\begin{array}{ccc}
0 &L^*(1) & 0\\ 0 & 0 &-L^*(1)^t \\ 0 & 0 & 0 \\
\end{array}\right),
\end{array}$$
where all $X,Y\in\Real^n$. We assume that each $\bar R_{xi}$ is zero on vectors
for which $\bar R_{xi}$ was not defined.

Since any holonomy algebra $\hol^\h$ of any other type with the orthogonal part $\h$ is contained
in $\hol^{1,\h}$, for the curvature  tensor $\bar R_x\in\R(\hol^\h)$  we have $\bar R_x\in\R(\hol^{1,\h})$,
and the decomposition for $\bar R_x$  can be obtained from the above decomposition and the  condition that  $\bar R_x$ 
takes values in $\hol^\h$.

For $\bar R_x\in\R(\hol^{2,\h})$ we have  
\begin{equation}
\bar R_{x4}=\bar R_{x5}=0 \text{ and } \bar R_{x}=\bar R_{x1}+\bar R_{x2}+\bar R_{x3}.\label{S42}
\end{equation}

For $\bar R_x\in\R(\hol^{3,\h,\varphi})$ we have  
\begin{equation}
\bar R_{x4}=0, \bar R_{x1}(X,Y)\in\ker\varphi  \text{ and } 
\bar R_{x5}(V,X)U=\varphi(\bar R_{x2}(V,X))U\text{  for all } X,Y\in\Real^n .\label{S43}
\end{equation}

For $\bar R_x\in\R(\hol^{4,\h,\psi})$ we have  
\begin{equation}
\bar R_{x4}=\bar R_{x5}=0, \bar R_{x1}(Z_1,Z_2)\in\ker\psi, \bar R_{x3}(V,Y)|_{\Real^{n_2}}=0,  \text{ and } 
\bar R_{x3}(V,X)Y=\eta(\psi(\bar R_{x2}(V,X)),Y)U\label{S44}\end{equation}
for all  $Z_1,Z_2\in\Real^n$ and $X+Y\in\Real^{n_1}\oplus\Real^{n_2}=\Real^{n}$.

\section{Decomposition of the curvature tensor on a manifold}\label{R}

In this section we use the notations of section 2. 
Let $x\in\M$.  For $T_x\M$ we have the decomposition $$T_x\M=\T_x\oplus S_x\oplus Tr_x.$$
In section \ref{G} we saw that having such decomposition of $T_x\M$, we can decompose the curvature tensor of $\M$  at the
point $x$ into 5 components, $\bar R_{x}=\bar R_{x1}+\bar R_{x2}+\bar R_{x3}+\bar R_{x4}+\bar R_{x5}$. Thus we have 
\begin{equation} \bar R=\bar R_1+\bar R_2+\bar R_3+\bar R_4+\bar R_5. \end{equation}
From the results of section \ref{G} it follows that for all $U,U_1\in \G(T^\bot)$, $X,Y,Z\in\G(S)$ and $V,V_1\in \G(Tr)$  we  have
\begin{align}
\bar R(X,Y)Z&=\underbrace{\bar R_1(X,Y)Z}_{\in \G(S)}+\underbrace{\bar R_2(X,Y)Z}_{\in \G(\T)},\label{l20}\\
\bar R(X,Y)U&=0,\label{5)}\\
\bar R(X,V)Y&=\underbrace{\bar R_2(X,V)Y}_{\in \G(S)}+\underbrace{\bar R_3(X,V)Y}_{\in \G(\T)},\label{l22}\\ 
\bar R(X,V)U&=\bar R_5(X,V)U\in \G(T^\bot), \label{1)}\\
\bar R(U,V)U_1&=\bar R_4(U,V)U_1\in \G(T^\bot), \label{4)}\\
\bar R(U,V)V_1&=\underbrace{\bar R_5(U,V)V_1}_{\in \G(S)}+\underbrace{\bar R_4(U,V)V_1}_{\in \G(Tr)},\label{3)}
\end{align}

Using this we will express $\bar R_1$, $\bar R_2$, $\bar R_3$, $\bar R_4$ and $\bar R_5$ in terms of the operators that we defined above.

Let $X,Y,Z\in \G(S)$, then  
\begin{align*}
\bar R(X,Y)Z&=\nb_X\nb_YZ-\nb_Y\nb_XZ-\nb_{[X,Y]}Z\overset{\eqref{l14}}{=}\n_X\n_YZ-\n_Y\n_XZ-\n_{[X,Y]}Z\\
&\overset{\eqref{l16}}{=}\n_X(\ns_YZ+\hs(Y,Z))-\n_Y(\ns_XZ-\hs(X,Z))-\ns_{[X,Y]}Z-\hs([X,Y],Z)\\
&\overset{(\ref{l16},\ref{l17})}{=}\ns_X\ns_YZ+ \hs(X,\ns_YZ)+\nst_X(\hs(Y,Z))\\&\qquad-\ns_Y\ns_XZ- \hs(Y,\ns_XZ)-\nst_Y(\hs(X,Z))
-\ns_{[X,Y]}Z-\hs([X,Y],Z).
\end{align*}
We have \begin{equation}\label{l25} [X,Y]=\nb_XY-\nb_YX=\n_XY-\n_YX=\ns_XY+\hs(X,Y) -\ns_YX-\hs(Y,X).\end{equation}
Thus, \begin{equation}\label{l26} \bar R(X,Y)Z=\underbrace{\Rs(X,Y)Z}_{\in \G(S)}+
\underbrace{(\ns_X\hs)(Y,Z)-(\ns_Y\hs)(X,Z)-\hs(\hs(X,Y)-\hs(Y,X),Z)}_{\in \G(\T)},\end{equation}
where $\Rs$ is the curvature tensor of the connection $\ns$ and
\begin{equation} (\ns_Xh)(Y,Z)=\nst_X(\hs(Y,Z))-\hs(\ns_XY,Z)-\hs(Y,\ns_XZ).\end{equation}

From \eqref{l20}  it follows that 
\begin{equation}\label{l27} \bar R_1(X,Y)Z=\Rs(X,Y)Z \end{equation}
and
\begin{equation}\label{l28} \bar R_2(X,Y)Z=(\ns_X\hs)(Y,Z)-(\ns_Y\hs)(X,Z)-\hs(\hs(X,Y)-\hs(Y,X),Z) \end{equation}
for all $X,Y,Z\in\G(S)$.

Note that from \eqref{l25} it follows that {\it the distribution $S$ is involutive if and only if}
$$\hs(X,Y)=\hs(Y,X) \text{ for all } X,Y\in\G(S).$$
From \eqref{l18} we see that this is equivalent to
$$\bar g(A_VX,Y)=\bar g(X,A_VY) \text{ for all } X,Y\in\G(S), V\in\G(Tr).$$

For $X,Y\in \G(S)$ and $V\in\G(Tr)$ we have   
\begin{align*}
\bar R(X,V)Y&=\nb_X\nb_VY-\nb_V\nb_XY-\nb_{[X,V]}Y\overset{\eqref{l14}}{=}\n_X\n_VY-\n_V\n_XY-\n_{[X,V]}Y\\
&\overset{\eqref{l16}}{=}\n_X(\ns_VY+\hs(V,Y))-\n_V(\ns_XY-\hs(X,Y))-\ns_{[X,V]}Y-\hs([X,V],Y)\\
&\overset{(\ref{l16},\ref{l17})}{=}\ns_X\ns_VY+ \hs(X,\ns_VY)+\nst_X(\hs(V,Y))\\&\qquad-\ns_V\ns_XY- \hs(V,\ns_XY)-\nst_V(\hs(X,Y))
-\ns_{[X,V]}Y-\hs([X,V],Y).
\end{align*}
We have $$[X,V]=\nb_XV-\nb_VX\overset{(\ref{l14},\ref{l15})}{=}-A_VX+\nt_XV-\ns_VX-\hs(V,X).$$
Thus, \begin{equation}\label{l29} \bar R(X,V)Y=\underbrace{\Rs(X,V)Y}_{\in \G(S)}+
\underbrace{(\n_X\hs)(V,Y)-(\ns_V\hs)(X,Y)+\hs(A_V,X)+\hs(\hs(V,X),Y)}_{\in \G(\T)},\end{equation}
where \begin{equation} (\n_X\hs)(V,Y)=\nst_X(\hs(V,Y))-\hs(\nt_XV,Y)-\hs(V,\ns_XY).\end{equation}

From  \eqref{l22} and  \eqref{l29}   it follows that 
\begin{equation}\label{l30} \bar R_2(X,V)Y=\Rs(X,V)Y \end{equation}
and
\begin{equation}\label{l31} \bar R_3(X,V)Y=(\n_X\hs)(V,Y)-(\ns_V\hs)(X,Y)+\hs(A_V,X)+\hs(\hs(V,X),Y)\end{equation}
for all $X,Y\in\G(S)$ and $V\in\G(Tr)$.

From \eqref{l27} and \eqref{l30} it follows that the curvature tensor $R^*$ can be found as follows 
\begin{equation}\label{l32}  \Rs(W_1,W_2)X=\pr_S(\bar R_1(W_1,W_2)X+\bar R_2(W_1,W_2)X) \text{ for all } W_1,W_2\in\G(T\M), X\in\G(S), \end{equation}
where $\pr_S$ is the projection on $S$ with respect to  decomposition \eqref{l13}. 

Let $x\in\M$ and  $M_x$ the lightlike hypersurface through $x$.
From \eqref{l32} it follows that for the curvature tensor of connection $\ns$ on the vector  bundle
$S(TM_x)=S|_{M_x}$ we have
\begin{equation}\label{l32B}
\Rs_{M_x}=\bar {R_1}|_{\G(TM_x)\times\G(TM_x)\times\G(S(TM_x))}. \end{equation}

From \eqref{l11} it follows that for the curvature tensor of $M_x$ we have
\begin{equation}\label{l32A} 
R_{M_x}=(\bar {R_1}+\bar {R_2})|_{\G(TM_x)\times\G(TM_x)\times\G(TM_x)}. \end{equation}
This shows that the curvature tensors $\bar R_1$  and $\bar R_2$  at a point $x\in\M$ depend only on
the lightlike  hypersurface through $x$, while the curvature  tensor  $\bar R_3$
depends on the links between different hypersurfaces.

For $U\in\G(\T)$ and $V,V_1\in\G(Tr)$ we have
\begin{align*}
\bar R(U,V)V_1&=\nb_U\nb_VV_1-\nb_V\nb_UV_1-\nb_{[U,V]}V_1\\
&\overset{\eqref{l15}}{=}\nb_U(-A_{V_1}V+\nt_VV_1)-\nb_V(-A_{V_1}U+\nt_UV_1)+A_{V_1}([U,V])-\nt_{[U,V]}V_1\\
&\overset{(\ref{l14},\ref{l15})}{=}-\n_UA_{V_1}V-A_{\nt_VV_1}U+\nt_U\nt_VV_1+\n_VA_{V_1}U+A_{\nt_UV_1}V\\
&\qquad -\nt_V\nt_UV_1+A_{V_1}(\nb_UV-\nb_VU)-\nt_{[U,V]}V_1\\
&\overset{(\ref{l16},\ref{l17})}{=}R^t(U,V)V_1-\ns_UA_{V_1}V-\hs(U,A_{V_1}V)-  A_{\nt_VV_1}U+\ns_VA_{V_1}U+\hs(V,A_{V_1}U)\\
&\qquad +A_{\nt_UV_1}V +A_{V_1}\nt_UV-A_{V_1}\nst_VU\\
&=\underbrace{R^t(U,V)V_1}_{\in\G(Tr)}-\underbrace{(\n_UA)_{V_1}V+(\n_VA)_{V_1}U}_{\in\G(S)}-\underbrace{\hs(U,A_{V_1}V)+\hs(V,A_{V_1}U)}_{\in\G(T^\bot)},
\end{align*}
where $$(\n_UA)_{V_1}V=\ns_UA_{V_1}V-A_{\nt_UV_1}V-A_{V_1}\nt_UV$$
and  $$(\n_VA)_{V_1}U=\ns_VA_{V_1}U-A_{\nt_VV_1}U-A_{V_1}\nst_VU.$$
From \eqref{3)} it follows that
\begin{equation}\label{R5UV} \bar R_5(U,V)V_1=(\n_VA)_{V_1}U-(\n_UA)_{V_1}V.\end{equation}
Using \eqref{4)} and \eqref{1)} we obtain
\begin{equation}\label{Rs1} \bar R_4(U,V)U_1=\Rs^t(U,V)U_1 \text{ and } \bar R_5(X,V)U=\Rs^t(X,V)U.\end{equation}
From \eqref{5)} it follows that \begin{equation}\label{Rs2} \Rs^t(X,Y)=0. \end{equation}
From \eqref{Rs1} and \eqref{Rs2} we obtain 
\begin{equation}\label{Rst} \Rs^t(W_1,W_2)U=\bar R_4(W_1,W_2)U+\bar R_5(W_1,W_2)U\end{equation}
for all $W_1,W_2\in\G(T\M)$ and $U\in\G(\T)$.

\section{Types of holonomy algebras}\label{type}
Let $\hol_x$ be the holonomy algebra of $\M$ at a point $x\in\M$. We suppose that $\hol_x$ is weakly irreducible.
In this section we give a criterion how to find the type of the holonomy algebra in terms of our global operators.

From section \ref{hol}, equations \eqref{Rst} and \eqref{R5UV} we obtain the following

\begin{prop}\label{prop1} The following conditions are equivalent

1) $\hol_x$ is of type 2 or 4;

2) $\bar R_4=\bar R_5=0$;

3) $\Rs^t=0$, i.e. the connection $\nst$  is flat;

4) $\bar R_4=0$ and $(\n_VA)_{V_1}U-(\n_UA)_{V_1}V=0$ for all
$U\in\G(\T)$, $V,V_1\in\G(Tr)$.
\end{prop} 

In the following proposition we use the vector bundle $\Hom(\so(S),\Real)$ over $\M$ such that $\Hom(\so(S),\Real)_y=\Hom(\so(S_y),\Real)$
for all $y\in\M$.  For a curve $\gamma$ in $\M$ we will denote by $\tau(\gamma)$ the parallel transport along $\gamma$. 

\begin{prop}\label{prop2}
The holonomy algebra $\hol_x$ is of type 3 if and only if the following conditions hold

\begin{itemize}
\item[1.] $\Rs^t(U,V)=0$ for all $U\in\G(T^\bot)$ and $V\in\G(Tr)$.
\item[2.] There exists a section $\varphi\in\G(\Hom(\so(S),\Real))$ such that 
\begin{itemize} 
\item[2.1.] $\Rs(X_y,Y_y)\in\ker\varphi_y$ for all $y\in\bar M$ and $X_y,Y_y\in S_y.$
\item[2.2.] $\Rs^t(V,X)=\varphi(\Rs(V,X))\cdot\id_{\G(T^\bot)}$ for all  $X\in \G(S)$ and $V\in\G(Tr)$.
\item[2.3.] There exist $y\in\M$, $X_y\in S_y$ and $V_y\in Tr_y$ such that $\Rs_y(V_y,X_y)\not\in\ker\varphi_y$.
\item[2.4.] For any  $y\in\M$ and any curve  $\gamma:[a,b]\to\M$ with $\gamma(a)=x$ and $\gamma(b)=y$ we have
$$\varphi_y(\Rs_y(V_y,X_y))=\varphi_x(\pr_{S_x}\circ\tau(\gamma)^{-1}\circ\Rs_y(V_y,X_y)\circ\tau(\gamma)|_{S_x})\text{ for all }
X_y\in S_y,V_y\in Tr_y.$$
\end{itemize}
\end{itemize} 
\end{prop}

{\it Proof.} Suppose that $\hol_x$ is of type 3. Since the holonomy algebras at different points of $\M$ are isomorphic, for any $y\in\M$
we have $\hol_y=\hol^{3,\h_y,\varphi_y}$, where $\h_y\subset\so(S_y)$ and $\varphi_y:\h_y\to \Real$ is a linear map. From this, \eqref{S43} and
\eqref{Rs1} follows the first statement of the proposition. We have $\so(S_y)=\h_y\oplus\h_y^\bot$, where $\h_y^\bot$ is the orthogonal complement
to $\h_y$ with respect to the Cartan-Killing form on $\so(S_y)$. Extend $\varphi_y$ to $\varphi_y:\so(S_y)\to\Real$ 
by setting $\varphi_y|_{\h_y^\bot}=0$. This gives a section $\varphi\in\G(\Hom(\so(S),\Real))$. Statements 2.1 and 2.2 follow from \eqref{S43},
\eqref{l27}, \eqref{l30} and \eqref{Rs1}. Suppose that 2.3 does not hold, then from statement 1, \eqref{Rs1} and proposition \ref{prop1}, it follows that
the holonomy algebra of $\M$ is of type 2 or 4, i.e. we obtain a contradiction.    
Let $y\in\M$ and let  $\gamma:[a,b]\to\M$ be a curve with $\gamma(a)=x$ and $\gamma(b)=y$.
For any $W_1$, $W_2\in T_x\M$ consider the operator
$$R^\gamma(W_1,W_2)=\tau(\gamma)^{-1}\circ\bar R_y(\tau(\gamma)W_1,\tau(\gamma)W_2))\circ\tau(\gamma):T_x\M\to T_x\M.$$
From the theorem of Ambrose and Singer it follows that $R^\gamma\in\R(\hol_x)$. Let $X_y\in S_y$, $V_y\in Tr_y$, 
and let $W_1=\tau(\gamma)^{-1}V_y\in T_x\M$, $W_2=\tau(\gamma)^{-1}X_y\in T_x$ (since $S\subset T$ and $T$ is parallel).
We can decompose $W_1$ and $W_2$ as follows $W_1=\pr_{T^\bot_x}W_1+\pr_{S_x}W_1+\pr_{Tr_x}W_1$ and $W_2=\pr_{T^\bot_x}W_2+\pr_{S_x}W_2$.
From section \ref{G} it follows that 
\begin{multline}\label{P21} R^\gamma_4(W_1,W_2)|_{T_x^\bot}=R^\gamma_4(\pr_{Tr_x}W_1,\pr_{S_x}W_2)|_{T_x^\bot}\\=
\varphi_x(R^\gamma_2(\pr_{Tr_x}W_1,\pr_{S_x}W_2))\cdot\id_{T_x^\bot}=
\varphi_x(\pr_{S_x}\circ R^\gamma_2(W_1,W_2)|_{S_x})\cdot\id_{T_x^\bot}.\end{multline}
From section \ref{G} it follows that the components  of any abstract curvature tensor   $R\in\R(\hol_x)$ can be obtained by
applying to $R$ some restrictions and projections.
We have \begin{multline}\label{P22}\varphi_x(\pr_{S_x}\circ R^\gamma_2 (W_1,W_2)|_{S_x})=\varphi_x(\pr_{S_x}\circ R^\gamma (W_1,W_2)|_{S_x})\\
=\varphi_x(\pr_{S_x}\circ\tau(\gamma)^{-1}\circ \bar R_y(V_y,X_y))\circ\tau(\gamma)|_{S_x})=
\varphi_x(\pr_{S_x}\circ\tau(\gamma)^{-1}\circ \bar R_{y2}(V_y,X_y))\circ\tau(\gamma)|_{S_x})\end{multline}
and \begin{multline*} R^\gamma_4(W_1,W_2)|_{T_x^\bot} = R^\gamma(W_1,W_2)|_{T_x^\bot}\\
=\tau(\gamma)^{-1}\circ \bar R_y(V_y,X_y)\circ\tau(\gamma)|_{T_x^\bot}=
\tau(\gamma)^{-1}\circ \bar R_{y4}(V_y,X_y)\circ\tau(\gamma)|_{T_x^\bot}.\end{multline*}  
Since $\tau(\gamma)T_x^\bot\subset T_y^\bot$, the last equality shows that 
$R^\gamma_4(W_1,W_2)|_{T_x^\bot}$ and $\bar R_{y4}(V_y,X_y)|_{T_y^\bot}$  
act on $T_x^\bot$ and $T_y^\bot$, respectively, as the multiplication on the same real number.
Since $\bar R_{y4}(V_y,X_y)|_{T_y^\bot}=\varphi_y(\bar R_{y2}(V_y,X_y))\cdot\id_{T_y^\bot}$, we see that
\begin{equation}\label{P23} R^\gamma_4(W_1,W_2)|_{T_x^\bot}=\varphi_y(\bar R_{y2}(V_y,X_y))\cdot\id_{T_x^\bot}.\end{equation}
Now statement 2.3 follows from \eqref{P21}, \eqref{P22}, \eqref{P23} and  \eqref{l30}.

Let us prove the  inverse. From \eqref{l30}, \eqref{Rs1} and proposition \ref{prop1} it follows that $\hol_x$
is of type 1 or 3. By the Ambrose and Singer theorem, the vector space $\hol_x$ is spanned
by the  elements $R^\gamma(W_1,W_2)$, where $W_1,W_2\in T_x\M$ and $\gamma$ is a curve in $\M$ with
the beginning at the point $x$. To show that $\hol_x$ is not of type 1, we must prove the
claim that if for some natural number $k$,  $\alpha_i\in\Real$, $W_{1i}$, $W_{2i}\in T_x\M$ and curves $\gamma_i$, where
$1\leq i\leq k$, holds $\sum_{i=1}^{k}\alpha_i R^{\gamma_i}(W_{1i},W_{2i})|_{T_x^\bot}\neq 0$, then
$\sum_{i=1}^{k}\alpha_i\pr_{S_x}\circ R^{\gamma_i}(W_{1i},W_{2i})|_{S_x}\neq 0$.
From section \ref{G} and statement 1 it follows that if for some $i$ holds $R^{\gamma_i}(W_{1i},W_{2i})|_{T_x^\bot}\neq 0$,
then $W_{1i}=V_i\in Tr_x$ and $W_{2i}=X_i\in S_x$ (or vice verse). 
As above we can show that from statement 2.4 it follows that $$R^{\gamma_i}(V_{i},X_{i})|_{T_x^\bot}=
\varphi_x(\pr_{S_x}\circ R^{\gamma_i}(W_{1i},W_{2i})|_{S_x}).$$ This proves the claim and the proposition. $\Box$

For any two sub-distributions $S_1,S_2\subset S$ denote by $\Hom(\so(S_1),S_2)$ the vector bundle over $\M$ 
such that $\Hom(\so(S_1),S_2)_y=\Hom(\so(S_{1y}),S_{2y})$ for all $y\in\M$. 

\begin{prop}\label{prop3}
The holonomy algebra $\hol_x$ is of type 4 if and only if 
it is of type 2 or of type 4 (see proposition \ref{prop1}) and
the following conditions hold

\begin{itemize}
\item[1.] There exist two parallel  sub-distributions $S_1,S_2\subset S$ such that $S=S_1\overset{\bot}{\oplus}S_2$ and the induced connection in $S_2$
is flat.
\item[2.] There exists a section $\psi\in\G(\Hom(\so(S_1),S_2))$ such that 
\begin{itemize} 
\item[2.1.] $\Rs(X_y,Y_y)\in\ker\psi_y$ for all $y\in\bar M$ and $X_y,Y_y\in S_y.$
\item[2.2.] $\bar R_3(V,Y)|_{S_2}=0$ and  $\bar R_3(V,X)Y=\bar g(\psi(\Rs(V,X)),Y)U$ for all  
$X\in \G(S_1)$, $Y\in\G(S_2)$, $U\in\G(T^\bot)$ and $V\in\G(Tr)$ such that $\bar g(U,V)=1$.
\item[2.3.] For any  $y\in\M$ and any curve  $\gamma:[a,b]\to\M$ with $\gamma(a)=x$ and $\gamma(b)=y$ we have
$$\bar g_x(\psi_y(\Rs_y(V_y,X_y)),\tau(\gamma)Y_x)
=\bar g_y(\psi_x(\pr_{S_{1x}}\circ\tau(\gamma)^{-1}\circ\Rs_y(V_y,X_y)\circ\tau(\gamma)|_{S_{1x}}),Y_x) $$  for all 
$X_y\in S_{1y}$, $Y_y\in S_{2y}$, and  $V_y\in Tr_y$.
\end{itemize}
\end{itemize} 
\end{prop}

The proof of proposition \ref{prop3} is similar to the proof of proposition \ref{prop2}.

Note that the statement 1 is equivalent to $$\Rs(W_1,W_2)\Gamma(S_1)\subset\Gamma(S_1) \text{ and } \Rs(W_1,W_2)|_{\G(S_2)}=0 \text{ for all } 
W_1,W_2\in\G(T_x\M).$$
Using \eqref{l31}, we can rewrite statement 2.2 in terms of $A$, $\hs$ and $\ns$. 

\vskip1cm

{\it Example 1.}
Consider the case when $(\M,\bar g)$ has an Abelian holonomy. It is known (see \cite{Sch}) that 
locally there exist coordinates $x^0,x^1,...,x^n,x^{n+1}$  
such that the metric $\bar g$ has the form $$\bar g=2dx^0dx^{n+1}+\sum_{i=1}^{n}(dx^i)^2+f\cdot(dx^{n+1})^2,$$
where $f$ is a function of $x^1,...,x^n,x^{n+1}$.
The nonzero Christoffel symbols are the following 
\begin{align*}
 \G^0_{n+1\, n+1 }&=\frac{1}{2}\frac{\partial f}{\partial x^{n+1}}, \\
 \G^i_{n+1\, n+1 }&=-\frac{1}{2}\frac{\partial f}{\partial x^{i}}, \\
\G^0_{i\, n+1 }&=\frac{1}{2}\frac{\partial f}{\,\, \partial x^{i}},
\end{align*}
where $1\leq i\leq n$. 

We see that for all $x$ from this coordinate neighborhood we have $\T_x=\Real(\frac{\partial}{\partial x^0})_x$ and 
we choose $S_x=\spa((\frac{\partial}{\partial x^1})_x,...,(\frac{\partial}{\partial x^n})_x)$. Then
$Tr_x=\Real N_x$, where $N=-\frac{1}{2}f \frac{\partial }{\partial x^{0}}+\frac{\partial}{\partial x^{n+1}}$.

Since for $1\leq i,j\leq n$ holds
$\n_{\frac{\partial}{\partial x^i}}\frac{\partial}{\partial x^j}=\n_{\frac{\partial}{\partial x^0}}\frac{\partial}{\partial x^j}=0$,
we have $$\ns_{\frac{\partial}{\partial x^i}}\frac{\partial}{\partial x^j}=\ns_{\frac{\partial}{\partial x^0}}\frac{\partial}{\partial x^j}=0 \text{ and } 
\hs\left(\frac{\partial}{\partial x^i},\frac{\partial}{\partial x^j}\right)=\hs\left(\frac{\partial}{\partial x^0},\frac{\partial}{\partial x^j}\right)=0.$$
For $1\leq j\leq n$ we  have $\n_N\frac{\partial}{\partial x^j}=\frac{1}{2}\frac{\partial f}{\partial x^{j}}\frac{\partial }{\partial x^{0}}$, hence
$$\ns_N\frac{\partial}{\partial x^j}=0 \text{ and } 
\hs\left(N,\frac{\partial}{\partial x^j}\right)=\frac{1}{2}\frac{\partial f}{\partial x^{j}}\frac{\partial }{\partial x^{0}}.$$

Furthermore, \begin{align*} \nb_NN&=\nb_{(-\frac{1}{2}f \frac{\partial }{\partial x^{0}}+\frac{\partial}{\partial x^{n+1}})}
(-\frac{1}{2}f \frac{\partial }{\partial x^{0}}+\frac{\partial}{\partial x^{n+1}})=
-\frac{1}{2}\frac{\partial f}{\partial x^{n+1}}\frac{\partial }{\partial x^{0}}+\frac{1}{2}\frac{\partial f}{\partial x^{n+1}}\frac{\partial }{\partial x^{0}}-
\sum_{i=1}^n\frac{1}{2}\frac{\partial f}{\partial x^{i}}\frac{\partial }{\partial x^{i}}\\
&=-\sum_{i=1}^n\frac{1}{2}\frac{\partial f}{\partial x^{i}}\frac{\partial }{\partial x^{i}},\\    
\nb_\frac{\partial}{\partial x^j}N&=\nb_\frac{\partial}{\partial x^j}(-\frac{1}{2}f \frac{\partial }{\partial x^{0}}+\frac{\partial}{\partial x^{n+1}})=
-\frac{1}{2}\frac{\partial f}{\partial x^{i}}\frac{\partial }{\partial x^{0}}+\frac{1}{2}\frac{\partial f}{\partial x^{i}}\frac{\partial }{\partial x^{0}}=0
\end{align*}
 and $$
\nb_\frac{\partial}{\partial x^0}N=0.$$
Hence, $\nt_WN=0$ for all $W\in\G(T\M)$, 
$$A_N\frac{\partial}{\partial x^0}=A_N\frac{\partial}{\partial x^i}=0 \text { and } 
A_NN=\sum_{i=1}^n\frac{1}{2}\frac{\partial f}{\partial x^{i}}\frac{\partial }{\partial x^{i}}.$$

Let $x\in\M$ be a point of this coordinate neighborhood and $M_x$ be  the lightlike hypersurface through $x$. 
Then for the shape operator  and the second fundamental form of $M_x$ we have
$$A=0 \text{ and } \hs=0$$
respectively.

For the curvature tensors we have $$\bar R_1=\bar R_2=\Rs=0 \text{ and } 
R_3(\frac{\partial }{\partial x^{i}},N)\frac{\partial }{\partial x^{j}}=- 
\frac{1}{2}\frac{\partial^2 f}{\partial x^{i}\partial x^{j}}\frac{\partial}{\partial x^0}.$$

\bibliographystyle{unsrt}

\end{document}